\documentclass[11pt]{amsart}
\usepackage{amssymb,amscd}

\newtheorem{thm}{Theorem}

\theoremstyle{remark}

\theoremstyle{definition}

\newcommand{\A}{{\mathcal{A}}}
\renewcommand{\O}{{\mathcal{O}}}

\newcommand{\C}{{\mathbb{C}}}

\newcommand{\Z}{{\mathbb{Z}}}

\newcommand{\End}{\mathrm{End}}

\newcommand{\gp}{\mathfrak p}
\def\zmod#1{\,\,({\rm mod}\,\,#1)}

\begin{document}
\title{A refined counter-example to the support conjecture for abelian varieties}
\author{Michael Larsen}
\email{larsen@math.indiana.edu}
\address{Department of Mathematics\\
    Indiana University \\
    Bloomington, IN 47405\\
    U.S.A.}

\author{Ren\'e Schoof}
\email{schoof@science.uva.nl}
\address{Universit\`a di Roma ``Tor Vergata'' \\
  Dipartimento di Matematica  \\
  Via della Ricerca Scientifica \\
  I-00133 Roma\\
  ITALY}

\thanks{The first named author was partially supported by NSF grant  DMS-0100537.}

\begin{abstract}
If $A/K$ is an abelian variety over a number field and $P$ and $Q$ are rational points, the original support conjecture asserted that if the order of $Q$ (mod $\gp$) divides the order of $P$
(mod $\gp$) for almost all primes $\gp$ of $K$, then $Q$ is obtained from $P$ by applying an endomorphism of $A$.  This is now known to be untrue.  In this note we prove that it is not even true modulo the torsion of $A$.
\end{abstract}
\maketitle

Let $A$ be an abelian variety over a number field $K$ and let $P$ and $Q$ be $K$-rational points of $A$.   By inverting a suitable element in the ring of integers of $K$, one can always find
a Dedekind domain $\O$ with fraction field $K$ such that $A$ extends to an abelian scheme 
$\A$ over $\O$ and $P$ and $Q$ extend to $\O$-points of $\A$.  Therefore, one can
speak of reducing $P$ and $Q$ (mod $\gp$) for almost all (i.e., all but finitely many) primes $\gp$.
In \cite{C-RS}, C.~Corrales-Rodrig\'a\~nez and R.~Schoof proved that when $\dim A = 1$,
the condition
\begin{equation}
\label{e:implication}
n P \equiv 0\pmod{\gp} \quad\Rightarrow\quad n Q\equiv 0 \pmod{\gp}
\end{equation}
for all integers $n$ and almost all prime ideals $\gp$ implies 
\begin{equation}
\label{e:strong}
Q =fP,\qquad\hbox{for some $f\in  \End_K(A)$.}
\end{equation}
In \cite{Larsen}, M.~Larsen proved 
that (\ref{e:implication}) does not imply (\ref{e:strong}) for general
abelian varieties but that it does imply 
\begin{equation}
\label{e:weak}
kQ =fP,\qquad\hbox{for some $f\in  \End_K(A)$}
\end{equation}
and some positive integer $k$.   The counter-example presented to (\ref{e:strong})
actually satisfies something stronger than (\ref{e:weak}), namely
\begin{equation}
\label{e:medium}
Q =fP+T,\qquad\hbox{for some $f\in  \End_K(A)$}
\end{equation}
and some torsion point $T\in A(K)$.\par
An early draft of \cite{arxiv} (version 2) claimed that (\ref{e:implication}) in fact implies (\ref{e:medium}).
The proof given was incorrect, and the statement was removed from subsequent versions.  (Version 3  is essentially the same as the published version \cite{Larsen}, while version 4 corrects a series of misprints, in which $P$ was written for $Q$ and vice versa throughout several paragraphs of the proof of the main theorem.)  
\par In this note we present an example to show that (\ref{e:implication}) does not imply
(\ref{e:medium}).

\begin{thm}
There exists an abelian variety $A$ over a number field $K$ and points $P$ and $Q$
which satisfy (\ref{e:implication}) but not (\ref{e:medium}).

\end{thm}

\begin{proof}
Let $p$ be a prime. Let $K$ be a number field and let $E$ be an elliptic curve over~$K$  without complex multiplication that possesses a point $R\in E(K)$ of infinite order.  Suppose in addition that the $p$-torsion points of~$E$ are rational over~$K$ and let $R_1,R_2\in E(K)$ be two independent points of order~$p$. Consider the abelian surface  $A$ obtained by dividing $E\times E$ by the subgroup generated by the point $(R_1,R_2)$. Then $A$ is defined over~$K$.  

We describe the ring of $K$-endomorphisms of~$A$.
Let $E(\C)\cong \C/\Lambda$, and let $\lambda_1,\lambda_2\in p^{-1}\Lambda$ map to
$R_1$ and $R_2$ respectively.  Thus, if for certain integers $a$ and $b$ we have  $a\lambda_1+b\lambda_2\in\Lambda$, then necessarily~$a,b\in p\Z$.  Let 
$$M = \Z\left(\begin{matrix}\lambda_1\\ \lambda_2\end{matrix}\right)+
\Lambda^2 \subset p^{-1}\Lambda^2.$$
The complex torus $A(\C)$ is isomorphic to $\C^2/M$, and any endomorphism 
of $A(\C)$ is given by a complex  $2\times 2$ matrix
\begin{equation}
\label{e:matrix}
\left(\begin{matrix}a&b\\c&d\end{matrix}\right)\in M_2(\C)
\end{equation}
with
\begin{equation}
\label{e:lambda-squared}
a \Lambda,b\Lambda,c\Lambda,d\Lambda\subset p^{-1}\Lambda
\end{equation}
and %
\begin{equation}
\label{e:torsion-image}
a\lambda_1+b\lambda_2\in k\lambda_1+\Lambda,
\ c\lambda_1+d\lambda_2\in k\lambda_2+\Lambda
\end{equation}
for some $k\in\Z$.
As $E$ does not have complex multiplication, (\ref{e:lambda-squared}) implies $pa,pb,pc,pd\in\Z$.
Multiplying (\ref{e:torsion-image}) by $p$, we deduce that $a,b,c,d\in\Z$,
and then (\ref{e:torsion-image}) implies $a-k,b,c,d-k\in p\Z$.
Conversely, any matrix
(\ref{e:matrix}) whose entries satisfy $a-k,b,c,d-k\in p\Z$ for some $k\in\Z$, 
lies in $\End(A(\C))$ and therefore
in $\End_{\C}A$. Since the curve $E$ and the points $R_1$, $R_2$ are defined over~$K$, it lies therefore in~$\End_K A$. 

Let $P$ and $Q$ denote the images of the points $(R,0)$ and $(R,R)$ in $A(K)$ respectively. Suppose that $nQ\equiv 0\zmod{\gp}$ for some prime $\gp$ of good reduction and characteristic different from~$p$.
This means that $(nR,nR)$ is contained in the subgroup generated by $(R_1,R_2)$ in the group of points on $E\times E$ modulo~$\gp$. Since the characteristic of~$\gp$ is not~$p$, the torsion points $R_1$ and $R_2$ are \textit{distinct}
modulo~$\gp$. This implies that  $nR\equiv 0\zmod{\gp}$. It follows that $nP\equiv 0\zmod{\gp}$. Therefore condition (\ref{e:implication}) is satisfied. And of course,  so is the conclusion (\ref{e:weak}) of Larsen's Theorem with $k=p$ and $f\in\End_K(A)$ the endomorphism with matrix $\begin{pmatrix}0&p\\0&0\\ \end{pmatrix}$.

However, (\ref{e:medium}) does not hold because there is no endomorphism $g\in\End_K(A)$ for which $P=gQ$ plus a torsion point.
Indeed, this would imply
$$
\begin{pmatrix}R\\0\\ \end{pmatrix}=\left[\begin{pmatrix}k&0\\0&k\\ \end{pmatrix}+p\begin{pmatrix}a&b\\ c&d\\ \end{pmatrix}\right]
\begin{pmatrix} R\\ R\\ \end{pmatrix} +\begin{pmatrix}T_1\\ T_2\\ \end{pmatrix}
$$
for some $k\in\Z$ and some torsion points $T_1,T_2\in E(K)$. Since $R$ has infinite order, inspection of the second coordinate shows that $k+pc+pd=0$ so that $k\equiv0\zmod p$. On the other hand,
looking at the first coordinate we see that $1=pa+pb+k$, a contradiction.
\end{proof}

\end{document}